\documentclass{amsart}

\usepackage{amssymb}
\usepackage{bezier}
\usepackage{epic,eepic}
\usepackage{epsfig} 
\usepackage{times}

\date{December 11, 2005}

\numberwithin{equation}{section} 
\newcommand{\la}{\langle} 
\newcommand{\ra}{\rangle} 

\newcommand{\R}{{\mathbb R}} 
\renewcommand{\phi}{\varphi} 
\renewcommand{\emptyset}{\varnothing} 
\newcommand{\cn}{\colon} 

\catcode`@=11 \@mparswitchfalse 

\newcounter{mnotecount}[section]

\swapnumbers

\newtheorem{thm}{Theorem}[section]
\newtheorem{lemma}[thm]{Lemma}
\newtheorem{prop}[thm]{Proposition}

\theoremstyle{definition}

\theoremstyle{remark}

\newcommand{\f}{\partial}

\renewcommand{\setminus}{\smallsetminus}
\newcommand{\s}{{\mathbb S}}

\newcommand{\hau}{{\mathcal H}}

\newcommand{\N}{\mathbb{N}}

\title{Nakajima's problem: convex bodies of constant 
  width and constant brightness}

\author{Ralph Howard}
\address{Department of Mathematics,
University of South Carolina,
Columbia, S.C. 29208, USA}
\email{howard\char'100math.sc.edu}
\urladdr{http://www.math.sc.edu/$\sim$howard}

\author{Daniel Hug}
\address{Mathematisches Institut, 
Universit{\"a}t Freiburg,
D-79104 Freiburg, 
Germany}
\email{daniel.hug@math.uni-freiburg.de}
\urladdr{http://home.mathematik.uni-freiburg.de/hug/}

\begin{document}

\dedicatory{Dedicated to Rolf Schneider on the occasion of his 65th birthday}

\begin{abstract}
For a convex body $K\subset\R^n$, the $k$th projection function of $K$
assigns to any $k$-dimensional linear subspace of $\R^n$ the
$k$-volume of the orthogonal projection of $K$ to that subspace.  Let
$K$ and $K_0$ be convex bodies in $\R^n$, and let $K_0$ be
centrally symmetric and satisfy a weak regularity and curvature
condition (which includes all $K_0$ with $\f K_0$ of class $C^2$ with
positive radii of curvature). Assume that $K$ and $K_0$
have proportional $1$st projection functions (i.e., width functions) 
and proportional $k$th projection functions.  For $2\le k<(n+1)/2$ and
for $k=3, n=5$ we show that $K$ and $K_0$ are homothetic. In the
special case where $K_0$ is a Euclidean ball, we thus obtain
characterizations of Euclidean balls as convex bodies of constant
width and constant $k$-brightness.
\end{abstract}

\maketitle

\section{Introduction and statement of results}

Let $K$ be a convex body (a compact, convex set with nonempty interior) 
in $\R^n$, $n\ge 3$. Assume that, for any line, the length of the 
projection of $K$ to the line is independent of that line and, for 
any hyperplane, the volume of the projection of $K$ to the hyperplane is 
independent of that hyperplane. Must $K$ then be a Euclidean ball? 

In dimension three, this problem has become known as Nakajima's
problem \cite{Nak}; see \cite{Chakrel}, \cite{ChakGroe}, \cite{Croft},
\cite{Gardnerbook}, \cite{Goodey}, \cite{Heil}. It is easy to check
that the answer to it is in the affirmative if $K$ is a convex body in
$\R^3$ of class $C^2$.  For general convex 
bodies in $\R^3$, the problem is much more difficult and a solution
has only been found  recently. Let $\mathbb{G}(n,k)$ denote the
Grassmannian of $k$-dimensional linear subspaces of $\R^n$. A convex
body $K$ in $\R^n$ is said to have {\em constant $k$-brightness},
$k\in\{1,\ldots,n-1\}$, if the $k$-volume $V_k(K\vert U)$ of the
orthogonal projection of $K$ to the linear subspace
$U\in\mathbb{G}(n,k)$ is independent of that subspace.  The map
$$
\pi_k\cn\mathbb{G}(n,k)\to\R,\qquad U\mapsto V_k(K\vert U),
$$ 
is referred to as the {\em $k$th projection function} of
$K$. Hence a convex body $K$ has constant width (i.e.\ constant
1-brightness) if it has constant $1$st projection function (width
function).

\begin{thm}[\cite{Howard:NP}]\label{ThmHoward}
Let $K$ be a convex body in $\R^n$ having constant width and constant
$2$-brightness. Then $K$ is a Euclidean ball.
\end{thm}

This theorem provides a complete solution of the 
Nakajima problem in $\R^3$ for general convex
bodies. In the present paper, we continue this line of research. Our
main result complements Theorem \ref{ThmHoward} by covering the cases
of convex bodies of constant width and constant $k$-brightness with 
$2\le k<(n+1)/2$ or $k=3$, $n=5$.
	
\begin{thm}\label{ThmHHnew}
Let $K$ be a convex body in $\R^n$ having  
constant width and constant $k$-brightness with  $2\le
k<(n+1)/2$, or $k=3, n=5$.
Then $K$ is a Euclidean ball. 
\end{thm}

The preceding two theorems can be generalized to pairs of convex
bodies $K,K_0$ having proportional projection functions, provided that
$K_0$ is centrally symmetric and has a minimal amount of regularity.

\begin{thm}\label{ThmHHnewgen} 
Let $K,K_0$ be convex bodies in $\R^n$, and let $K_0$ be centrally
symmetric with positive principal radii of curvature on some Borel
subset of the unit sphere of positive measure.  Let $2\le
k<(n+1)/2$, or let $k=3,n=5$ in which case assume  the surface area
measure $S_4(K_0,\cdot)$ of $K_0$ is absolutely continuous with positive density.  
Assume that there
are constants $\alpha,\beta>0$ such that
$$
\pi_1(K)=\alpha\, \pi_1(K_0)\qquad\text{and}\qquad \pi_k(K)=\beta\, \pi_k(K_0).
$$
Then $K$ and $K_0$ are homothetic.
\end{thm}

As the natural measure on the unit sphere, $\s^{n-1}$, we use the invariant
Haar probability measure (i.e.\ spherical Lebesgue measure), 
or what is the same thing the $(n-1)$-dimensional Hausdorff measure, $\mathcal{H}^{n-1}$, 
normalized so that the total mass is one. 
We view the principal radii of curvature 
as functions of the unit normal, despite the fact that the unit normal map is in 
general a set valued function (cf.\ the beginning of Section 2 below). The assumption that 
the principal radii of curvature are positive on a set of positive measure 
means that there is a Borel subset of $\s^{n-1}$ of positive measure such that 
on this set the 
reverse Gauss map is single valued, differentiable (in a generalized sense) and 
the eigenvalues of the differential are positive. Explicitly, 
this condition can be stated in terms of second order
differentiability properties of the support function (again see
Section~\ref{sec:prelim}).  In particular, it is certainly satisfied 
if $K_0$ is of
class $C^2_+$, and therefore letting $K_0$ be a Euclidean ball recovers
Theorem~\ref{ThmHHnew}. The required condition 
allows for parts of $K_0$ to be
quite irregular.  For example if $\f K_0$ has a point that has a small
neighborhood where $\f K_0$ is $C^2$ with positive Gauss-Kronecker
curvature, then the assumption will hold, regardless of how rough the
rest of the boundary is.  For example a ``spherical polyhedron''
constructed by intersecting a finite number of Euclidean balls in
$\R^n$ will satisfy the condition.  More generally if the convex body
$K_0$ is an intersection of a finite collection of bodies of class
$C^2_+$, it will satisfy the condition.

Theorem~\ref{ThmHHnewgen} extends the main results in
\cite{HH:NPsmooth} for the range of dimensions
$k,n$ where it applies by reducing the regularity assumption on $K_0$
and doing away with any regularity assumptions on $K$.  However, the
classical Nakajima problem, which concerns the case $n=3$ and $k=2$,
is not covered by the present approach.

Despite recent progress on the Nakajima problem various questions
remain open.  For instance, can Euclidean balls be characterized as
convex bodies having constant width and constant $(n-1)$-brightness if
$n\ge 4$? This question is apparently unresolved even for smooth
convex bodies. A positive answer is available for smooth convex bodies
of revolution (cf.~\cite{HH:NPsmooth}).  From the arguments of the
present paper the following proposition is easy to check.

\begin{prop}
Let $K,K_0\subset\R^n$ be convex bodies that have a common axis of
revolution. Let $K_0$ be centrally symmetric with positive principal 
radii of curvature almost everywhere. Assume that $K$ and $K_0$ have
proportional width functions and proportional $k$th projection
functions for some $k\in\{2,\ldots,n-2\}$.  Then $K$ and $K_0$ are
homothetic.
\end{prop} 

It is a pleasure for the authors to dedicate this paper to 
Rolf Schneider. Professor Rolf Schneider has been a large source 
of inspiration for countless students and colleagues all over the 
world. His willingness to communicate and share his 
knowledge make contact with him a pleasurable and mathematically 
rewarding experience.  The second named author has particularly 
been enjoying many years of support, personal interaction and joint research.

\section{Preliminaries}\label{sec:prelim}

Let $K$ be a convex body in $\R^n$, and let $h_K\cn \R^n\to \R$ be the
support function of $K$, which is a convex function.  For $x\in \R^n$
let $\f h_K(x)$ be the subdifferential of $h_K$ at $x$.  This is the
set of vectors $v\in \R^n$ such that the function $h_K-\la v,\cdot\ra$
achieves its minimum at $x$.  It is well known that, for all $x\in \R^n$, 
$\f h_K(x)$ is a nonempty compact convex set and is a singleton
precisely at those points where $h_K$ is differentiable in the
classical sense (cf.~\cite[pp.~30--31]{Schneider1993}).  For $u\in
\s^{n-1}$ the set $\f h_K(u)$ is exactly the set of $x\in \f K$ such
that $u$ is an outward pointing normal to $K$ at $x$
(cf.~\cite[Thm~1.7.4]{Schneider1993}).  But this is just the
definition of the reverse Gauss map (which in general is not single
valued, but a set valued function) and so the function $u\mapsto \f
h_K(u)$ gives a formula for the reverse Gauss map in terms of the
support function.

In the following,  by ``almost everywhere'' on the unit 
sphere or by ``for almost all unit vectors'' we mean for all unit 
vectors with the possible exclusion of 
a set of spherical Lebesgue measure zero. 
 A theorem of Aleksandrov states that a convex function has a
generalized second derivative almost everywhere, which we will view as
a positive semidefinite symmetric linear map rather than a symmetric
bilinear form.  This generalized derivative can either be defined in
terms of a second order approximating Taylor polynomial at the point,
or in terms of the set valued function $x\mapsto \f h_K(x)$ being
differentiable in the sense of set valued functions (both these
definitions are discussed in \cite[p.~32]{Schneider1993}).  At points
where the Aleksandrov second derivative exists $\f h_K$ is single
valued.  Because $h_K$ is positively homogeneous of degree one, if it
is Aleksandrov differentiable at a point $x$, then it is Aleksandrov
differentiable at all points $\lambda x$ with $\lambda>0$.  Then
Fubini's theorem implies that not only is $h_K$ Aleksandrov
differentiable at $\mathcal{H}^n$ almost all points of $\R^n$, but it is also
Aleksandrov differentiable at $\hau^{n-1}$ almost all points of
$\s^{n-1}$.  For points $u\in \s^{n-1}$ where it exists, let
$d^2h_K(u)$ denote the Aleksandrov second derivative of $h_K$. Let $u^\perp$ 
denote the orthogonal complement of $u$. Then
the restriction $d^2h_K(u) \vert u^\perp$ is the derivative of
the reverse Gauss map at $u$.  The eigenvalues of $ d^2h_K(u) \vert
u^\perp$ are the principal radii of curvature at $u$.  As the
discussion above shows these exist at almost all points of $\s^{n-1}$.

A useful tool for the study of projection functions of
convex bodies are the surface area measures. An introduction to these
Borel measures on the unit sphere is given in \cite{Schneider1993}, a
more specialized reference (for the present purpose) is contained in
the preceding work \cite{HH:NPsmooth}. The top order surface area
measure $S_{n-1}(K,\cdot)$ of the convex body $K\subset\R^n$ can be
obtained as the $(n-1)$-dimensional Hausdorff measure
$\mathcal{H}^{n-1}$ of the reverse spherical image of Borel sets of
the unit sphere $\mathbb{S}^{n-1}$.  The Radon-Nikodym derivative of
$S_{n-1}(K,\cdot)$ with respect to the spherical Lebesgue measure is
the product of the principal radii of curvature of $K$.  Since for almost 
every $u\in \s^{n-1}$, the
radii of curvature of $K$ at $u\in\mathbb{S}^{n-1}$ are the
eigenvalues of $d^2h_K(u)\vert u^\perp$, the Radon-Nikodym derivative
of $S_{n-1}(K,\cdot)$ with respect to spherical Lebesgue measure is
the function $u\mapsto \det\left( d^2h_K(u) \vert u^\perp\right)$, which 
is defined almost everywhere on $\s^{n-1}$. In particular, if $S_{n-1}(K,\cdot)$ 
is absolutely continuous with respect to spherical Lebesgue measure, the 
density function is just the Radon-Nikodym derivative. 
For explicit definitions of these and other basic notions of convex
geometry needed here, we refer to \cite{Schneider1993} and
\cite{HH:NPsmooth}.

The following lemma contains more precise information about the
Radon-Nikodym derivative of the top order surface area measure.  We
denote the support function of a convex body $K$ by $h$, if $K$ is
clear from the context. For a fixed unit vector $u\in\mathbb{S}^{n-1}$
and $i\in\N$, we also put $ \omega_i:=\left\{v\in
\mathbb{S}^{n-1}:\langle v,u\rangle\ge 1-(2i^2)^{-1}\right\}$,
whenever $u$ is clear from the context.  Hence $\omega_i\downarrow
\{u\}$, as $i\to\infty$, in the sense of Hausdorff convergence of
closed sets.

\begin{lemma}\label{RND}
Let $K\subset\R^n$ be a convex body. If $u\in\mathbb{S}^{n-1}$ is a point of second 
order differentiability of the support function $h$ of $K$, then
$$
\lim_{i\to\infty}\frac{S_{n-1}(K,\omega_i)}{\mathcal{H}^{n-1}(\omega_i)}=
\det\left(d^2h(u)\vert u^\perp\right).
$$
\end{lemma}

\begin{proof} This is implicitly contained in the proof of Hilfssatz 2 in 
\cite{Leichtweiss88}. A similar argument,  
in a slightly more involved situation, can be found in \cite{Hug.aff}.
\end{proof}

An analogue of Lemma \ref{RND} for curvature measures is provided in 
\cite[(3.6) Hilfssatz]{Schneider1979}.

\bigskip

As another ingredient in our approach to Nakajima's problem, we need
two simple algebraic lemmas.  Here we write $|M|$ for the cardinality
of a set $M$.  If $x_1,\dots,x_n$ are real numbers and
$I=\{i_1,\dots,i_k\}\subseteq \{1,\dots,n\}$ we set
$x_I:=x_{i_1}\dots x_{i_k}$. We also put $x_{\emptyset}:=1$.

\begin{lemma}\label{alg2}
Let $b>0$ be fixed. Let $x_1,\ldots,x_{n-1},y_1,\ldots,y_{n-1}$ be
nonnegative real numbers satisfying
$$
x_i+y_i=2\qquad \text{and}\qquad
x_I+y_I=2b
$$ for all $i=1,\ldots,n-1$ and all $I\subset\{1,\ldots,n-1\}$ with
$|I|=k$, where $k\in\{2,\ldots,n-2\}$.  Then
$|\{x_1,\ldots,x_{n-1}\}|\le 2$ and $|\{y_1,\ldots,y_{n-1}\}|\le 2$.
\end{lemma}

\begin{proof} We can assume that $x_1\le\dots\le x_{n-1}$. Then we
  have $y_1\ge \dots\ge y_{n-1}$.

\bigskip

If $x_1=0$, then $y_1=2$. Further, for $I'\subset\{2,\ldots,n-1\}$
with $|I'|=k-1$, we have $y_1y_{I'}=2b$, hence $y_{I'}=b$. Since $k\ge
2$, we get $y_2,\ldots,y_{n-1}>0$. Moreover, since $k-1\le n-3$, we
conclude that $y_2=\dots=y_{n-1}$. This shows that also
$x_2=\dots=x_{n-1}$, and thus $|\{x_1,\ldots,x_{n-1}\}|\le 2$ and
$|\{y_1,\ldots,y_{n-1}\}|\le 2$.

\bigskip

If $y_{n-1}=0$, the same conclusion is obtained by symmetry.

\bigskip

If $x_1>0$ and $y_{n-1}>0$, then
$x_1,\ldots,x_{n-1},y_1,\ldots,y_{n-1}>0$. Now we fix any set
$J\subseteq\{1,\ldots,n-1\}$ with $|J|=k+1$. The argument at the
beginning of the proof of Lemma 4.2 in \cite{HH:NPsmooth} shows that
$|\{x_i:i\in J\}|\le 2$. Since $k+1\ge 3$, we first obtain that
$|\{x_1,\ldots,x_{n-1}\}|\le 2$, and then also
$|\{y_1,\ldots,y_{n-1}\}|\le 2$.
\end{proof}

\begin{lemma}\label{alg3}
Let $n\ge 4$, and let $b>0$ be fixed. 
Let $x_1,\ldots,x_{n-1},y_1,\ldots,y_{n-1}$ be nonnegative real numbers satisfying 
$$
x_i+y_i=2\qquad \text{and}\qquad
x_I+y_I=2b
$$
for all $i=1,\ldots,n-1$ and all $I\subset\{1,\ldots,n-1\}$ with $|I|=n-2$. 
Then
\begin{equation}\label{newstar}
\prod_{l\neq i,j}x_l=\prod_{l\neq i,j}y_l=b
\end{equation}
whenever $i,j\in\{1,\ldots,n-1\}$ are such that $x_i\neq x_j$. 
\end{lemma}

\begin{proof}
For the proof, we may assume that $i=1$ and $j=n-1$, to simplify the
notation. Then we have
\begin{align*}
x_1\cdots x_{n-2}+y_1\cdots y_{n-2}&=2b,\\
x_2\cdots x_{n-1}+y_2\cdots y_{n-1}&=2b,
\end{align*}
which implies that
$$
x_2\cdots x_{n-2}(x_{n-1}-x_1)+y_2\cdots y_{n-2}(y_{n-1}-y_1)=0.
$$
Moreover, $x_1+y_1=2=x_{n-1}+y_{n-1}$ yields
$$
x_{n-1}-x_1=y_1-y_{n-1}\neq 0,
$$
and thus
$$
x_2\cdots x_{n-2}=y_2\cdots y_{n-2}.
$$
Hence
\begin{align*}
2x_2\cdots x_{n-2}&=(x_1+y_1)x_2\cdots x_{n-2}=x_1x_2\cdots x_{n-2}+y_1x_2\cdots x_{n-2}\\
&=x_1x_2\cdots x_{n-2}+y_1y_2\cdots y_{n-2}=2b,
\end{align*}
and thus
$$
b=x_2\cdots x_{n-2}=y_2\cdots y_{n-2}.
$$
\end{proof}

\section{Proofs}

First, by possibly dilating $K$, we can assume that $\alpha=1$.  Hence
the assumption can be stated as
\begin{equation}\label{a1}
\pi_1(K)= \pi_1(K_0)\qquad\text{and}\qquad \pi_k(K)=\beta\, \pi_k(K_0)
\end{equation}
for some $k\in\{2,\ldots,n-2\}$. Let $K^*$ denote the reflection of
$K$ in the origin.  Then \eqref{a1} yields that
$$
K+K^*=2K_0\qquad\text{and}\qquad
V_{k}(K\vert U)=\beta\, V_{k}(K_0\vert U)
$$ for all $U\in\mathbb{G}(n,k)$. Minkowski's inequality (cf.\
\cite{Schneider1993}) then implies that
\begin{align*}
V_{k}(2K_0\vert U)=\,&V_{k}(K\vert U+K^*\vert U)\\
\ge\,&\left(V_{k}(K\vert U)^{\frac{1}{k}}+V_{k}
(K^*\vert U)^{\frac{1}{k}}\right)^{k}\\
=\,&\left(2V_{k}(K\vert U)^{\frac{1}{k}}\right)^{k}\\
=\,&\beta\, V_{k}(2K_0\vert U).
\end{align*}
Equality in Minkowski's inequality will hold if and only if $K^*\vert
U$ and $K\vert U$ are homothetic.  As they have the same volume this
is equivalent to their being translates of each other, in which case
$K\vert U$ is centrally symmetric. 
Hence $\beta\le 1$ with equality if and only if $K\vert
U$ is centrally symmetric for all linear subspaces $U\in
\mathbb{G}(n,k)$. Since $k\ge 2$, this is the case if and only if $K$
is centrally symmetric
(cf.~\cite[Thm.~3.1.3]{Gardnerbook}). 
So if $\beta=1$, then $K$ and $K_0$ must be homothetic.

In the following, we assume that $\beta\in (0,1)$.  This will lead to
a contradiction and thus prove the theorem.

We write $h,h_0$ for the support functions of $K,K_0$. Here and in
the following, ``almost all'' or ``almost every'' refers to the
natural Haar probability measure on $\s^{n-1}$.  Moreover a linear
subspace ``$E$'' as an upper index indicates that the corresponding
functional or measure is considered with respect to $E$ as the
surrounding space.  By assumption there is a Borel subset $P\subseteq
\s^{n-1}$ with positive measure such that for all $u\in P$ all the
radii of curvature of $K_0$ in the direction $u$ exist and are
positive.  As $K_0$ is symmetric we can assume that $u\in P$ if and
only if $-u\in P$.  Let $N$ be the set of points $u\in \s^{n-1}$ where the
principal radii of curvature of $K$ do not exist. Since $N$ is
the set of points where the Alexandrov second derivative of $h$ does
not exist, it is a set of measure zero.  By replacing $P$ by
$P\setminus(N\cup (-N))$ we can assume that the radii of curvature of
both $K_0$ and $K$ exist at all points of $P$.  As both $N$ and $-N$
have measure zero this set will still have positive measure.

Let 
$u\in \mathbb{S}^{n-1}$ be such that $h$ and $h_0$ are second
order differentiable at $u$ and at $-u$ and that the radii of
curvature of $K_0$ at $u$ are positive.  This is true of all points
$u\in P$, which is not empty as it has positive measure.  Let $E\in
\mathbb{G}(n,k+1)$ be such that $u\in E$.  Then the assumption implies
that also
$$
\pi^E_k(K\vert E)=\beta\, \pi^E_k(K_0\vert E).
$$
Hence we conclude as in \cite{HH:NPsmooth} that 
$$
S_k^E(K\vert E,\cdot)+ S_k^E(K^*\vert E,\cdot)=2\beta\, S_k^E(K_0\vert E,\cdot).
$$ Since $h(K\vert E,\cdot)=h_K\vert E$ and $h(K_0\vert
E,\cdot)=h_{K_0}\vert E$ are second order differentiable at $u$ and at
$-u$ with respect to $E$, Lemma \ref{RND} applied with respect to the
subspace $E$ implies that
\begin{multline*}
\det\left(d^2h_{K\vert E}(u)\vert E\cap u^\perp\right)+
\det\left(d^2h_{K^*\vert E}(u)\vert E\cap u^\perp\right)\\
=2\beta \,
\det\left(d^2h_{K_0\vert E}(u)\vert E\cap u^\perp\right).
\end{multline*}
Since $h$ and $h_0$ are second order differentiable at $u$ and at
$-u$, the linear maps 
$$
L(h)(u)\cn T_u\mathbb{S}^{n-1}\to T_u\mathbb{S}^{n-1}, \quad v\mapsto d^2h(u)(v),
$$
$$
L(h_0)(u)\cn T_u\mathbb{S}^{n-1}\to T_u\mathbb{S}^{n-1}, \quad v\mapsto d^2h_0(u)(v),
$$ 
are well defined and positive semidefinite.  
Since the radii of curvature of $K_0$ at $u$ are
positive, we can define
$$
L_{h_0}(h)(u):=L(h_0)(u)^{-1/2}\circ L(h)(u)\circ L(h_0)(u)^{-1/2}
$$
as in \cite{HH:NPsmooth} in the smooth case.

In this situation, the arguments in \cite{HH:NPsmooth} can be repeated
to yield that
\begin{equation}\label{star2}
\begin{aligned}
L_{h_0}(h)(u) + L_{h_0}(h)(-u) =&\, 2\,{\rm id}\\
\wedge^{k}L_{h_0}(h)(u) + \wedge^{k}L_{h_0}(h)(-u) =&\, 2\beta\, \wedge^{k} {\rm id},
\end{aligned}
\end{equation}
where ${\rm id}$ is the identity map on $T_u\mathbb{S}^{n-1}$.  Lemma 3.4 in
\cite{HH:NPsmooth} shows that $L_{h_0}(h)(u)$ and $L_{h_0}(h)(-u)$
have a common orthonormal basis of eigenvectors $e_1,\ldots,e_{n-1}$,
with corresponding eigenvalues (relative principal radii of curvature)
$x_1,\ldots,x_{n-1}$ at $u$ and with eigenvalues  
$y_1,\ldots,y_{n-1}$ at $-u$.  After a
change of notation (if necessary), we can assume that $0\le x_1\le
x_2\le \dots\le x_{n-1}$.  By \eqref{star2} we thus obtain
\begin{equation}\label{3.2s}
x_i+y_i=2\qquad\text{and}\qquad x_I+y_I=2\beta
\end{equation}
for $ i=1,\ldots,n-1$ and $I\subset\{1,\ldots,n-1\}$ with $|I|=k$.

\bigskip

{\bf Proof of Theorem \ref{ThmHHnewgen} when $2\le k < 
(n+1)/2$.}\   From \eqref{3.2s} and Lemma~\ref{alg2} we conclude that
there is some $\ell\in\{0,\ldots,n-1\}$ such that
$$
x_1=\dots=x_{\ell}<x_{\ell+1}=\dots=x_{n-1}
\qquad\text{and}\qquad 
y_1=\dots=y_{\ell}>y_{\ell+1}=\dots=y_{n-1}.
$$

(a) If $k\le \ell$, then
$$
x_1+y_1=2\qquad\text{and}\qquad x_1^k+y_1^k=2\beta .
$$
Hence 
$$
1=\left(\frac{x_1+y_1}{2}\right)^k\le \frac{x_1^k+y_1^k}{2}=\beta,
$$
contradicting the assumption that $\beta<1$.

\medskip

(b) Let $k>\ell$. Since $k<(n+1)/2$ we have $2k<n+1$ or
$k<n+1-k$. Hence $k\le n-k<n-\ell$, and thus $k\le n-1-\ell$. But then
$$
x_{\ell+1}+y_{\ell +1}=2\qquad\text{and}\qquad x_{\ell +1}^k+y_{\ell +1}^k=2\beta ,
$$ and we arrive at a contradiction as before.  This proves Theorem
\ref{ThmHHnewgen} when $2\le k < (n+1)/2$\qed

\bigskip

{\bf Proof of Theorem \ref{ThmHHnewgen} when $k=3, n=5$.} 
In this case we are assuming that $K_0$ has positive radii of curvature at
almost all points of $\s^{n-1}$.  As $h$ has Alexandrov second
derivatives at almost all points, for almost all
$u\in\s^{n-1}$ the radii of curvature of $K$ exist at both $u$ and
$-u$ and at these unit vectors $K_0$ has positive radii of curvature. 
Recall that $x_1\le\dots\le x_4$ are the eigenvalues of
$L_{h_0}(h)(u)$.  We distinguish three cases each of which will lead to a
contradiction.

\bigskip

(a) $x_1\neq x_2$. Then Lemma \ref{alg2} yields that $x_1<x_2=x_3=x_4$ and 
therefore also $y_2=y_3=y_4$. Hence
$$
x_2^3+y_2^3=2\beta\qquad\text{and}\qquad x_2+y_2=2,
$$
and thus
$$
1=\left(\frac{x_2+y_2}{2}\right)^3\le \frac{x_2^3+y_2^3}{2}=\beta,
$$
contradicting that $\beta<1$.  So this case can not arise.

\bigskip

(b) $x_1=x_2$ and $x_1=x_3$, i.e.\ $x_1=x_2=x_3$. Then also $y_1=y_2=y_3$, and we get
$$
x_1^3+y_1^3=2\beta\qquad\text{and}\qquad x_1+y_1=2,
$$
which, as before, leads to a contradiction and thus this case can not
arise. 

\bigskip

(c) $x_1=x_2$ and $x_1\neq x_3$, i.e.\ $x_1=x_2<x_3= x_4$ by Lemma
\ref{alg2}. Since $x_1\neq x_3$, Lemma \ref{alg3} implies that
\begin{equation}\label{i}
x_2x_4=\beta=y_2y_4. 
\end{equation}
In addition, we have
\begin{equation}\label{ii}
x_2+y_2=2=x_4+y_4.
\end{equation}
We show that these equations determine $x_2,x_4,y_2,y_4$ as functions
of $\beta$. Substituting \eqref{i} into \eqref{ii}, we get
$$
\frac{\beta}{x_4}+y_2=2,\qquad x_4+\frac{\beta}{y_2}=2.
$$
Combining these two equations, we arrive at
$$
y_2+\frac{\beta}{2-\frac{\beta}{y_2}}=2,
$$
where we used that $x_4=2-\frac{\beta}{y_2}\neq 0$. This 
equation for $y_2$ can be rewritten as
$$
y_2^2-2y_2+\beta=0.
$$
Hence, we find that (recall that $0<\beta<1$)
$$
y_2=1\pm \sqrt{1-\beta}.
$$
Consequently,
$$
x_2=2-y_2=1\mp\sqrt{1-\beta}.
$$
From \eqref{i}, we also get
$$
x_4=\frac{\beta}{x_2}=\frac{\beta}{1\mp\sqrt{1-\beta}}=1\pm \sqrt{1-\beta},
$$
and finally again by \eqref{i}
$$
y_4=\frac{\beta}{y_2}=\frac{\beta}{1\pm\sqrt{1-\beta}}=1\mp \sqrt{1-\beta}.
$$
Since $x_1=x_2< x_3=x_4$, this shows that
\begin{equation}\label{radcurv}
x_1=x_2=1-\sqrt{1-\beta},\qquad x_3=x_4=1+\sqrt{1-\beta}.
\end{equation}
By assumption the surface area measure $S_{4}(K_0,\cdot)$ of $K_0$ is
absolutely continuous with density function $u\mapsto
\det(d^2h_0(u)\vert u^\perp)$.  Since $K+K^*=2K_0$, the non-negativity
of the mixed surface area measures
$S(K[i],K^*[4-i],\cdot)$ and the multilinearity
of the surface area measures yields that
\begin{align*}
S_{4}(K,\cdot)\le\,&\sum_{i=0}^{4}\binom{4}{i}S(K[i],K^*[4-i],\cdot)\\
=\,&S_{4}(K+K^*,\cdot)=2^{4}\, S_{4}(K_0,\cdot).
\end{align*}
This implies that $S_{4}(K,\cdot)$ is absolutely continuous as well,
with density function $ u\mapsto \det(d^2h(u)\vert u^\perp)$. Now
observe that the cases (a) and (b) have already been excluded and
therefore the present case (c) is the only remaining one. Hence, using
the definition of $L_{h_0}(h)(u)$,
$$
\frac{\det(d^2h(u)\vert u^\perp)}{\det(d^2h_0(u)\vert u^\perp)}
  =\det(L_{h_0}(h)(u))=x_1x_2x_3x_4 =\beta^2,
$$
for almost all  $u\in \mathbb{S}^{4}$. Thus we deduce that
$$
S_{4}(K,\cdot)=\beta^2\, S_{4}(K_0,\cdot).
$$
Minkowski's uniqueness theorem now implies that $K$ and $K_0$ are homothetic, 
hence $K$ is centrally symmetric.  Symmetric convex bodies with
the same width function are translates of each other.  But then again 
$\beta=1$, a contradiction.


\providecommand{\bysame}{\leavevmode\hbox to3em{\hrulefill}\thinspace}
\providecommand{\MR}{\relax\ifhmode\unskip\space\fi MR }
\providecommand{\MRhref}[2]{%
  \href{http://www.ams.org/mathscinet-getitem?mr=#1}{#2}
}
\providecommand{\href}[2]{#2}

\end{document}